\documentclass[11pt]{article}

\usepackage{amsmath,amsthm,amsfonts,amssymb,mathrsfs,amscd}
\theoremstyle{plain}
\newtheorem{Theorem}{Theorem}
\newtheorem{statement}{Statement}
\newtheorem{Lemma}{Lemma}
\newtheorem{corollary}{Corollary}

\newtheorem*{littlewood}{Littlewood conjecture}
\newtheorem*{oppenheim}{Oppenheim conjecture on linear forms}
\newtheorem*{sailed_oppenheim}{Reformulated Oppenheim conjecture}
\newtheorem*{sailed_oppenheim_3dim}{Reformulated three--dimensional Oppenheim conjecture}

\theoremstyle{definition}
\newtheorem{defin}{Definition}
\newtheorem{Example}{Example}

\theoremstyle{remark}
\newtheorem{Remark}{Remark}
\newtheorem{nproof}{}

\DeclareMathOperator{\SL}{SL}
\newcommand\oSt{\mathrm{St}}
\DeclareMathOperator{\Aff}{Aff}
\DeclareMathOperator{\conv}{conv}
\DeclareMathOperator{\opint}{int}
\DeclareMathOperator{\aff}{aff}
\DeclareMathOperator{\opvol}{vol}
\DeclareMathOperator{\opext}{ext}
\DeclareMathOperator{\GL}{GL}

\let\leq=\leqslant
\let\geq=\geqslant

\author{O.~N.~German, E.~L.~Lakshtanov}
\date{12.07.2006}

\title{On multidimensional generalization of the Lagrange theorem on continued fractions
\thanks{The research of the first author was supported financially by RFBR
(grant no.~06-01-00518), INTAS (grant no.~03-51-5070) and the Programme of the President of Russian
Federation (grant no.~MK-6370.2006.1), the research of the second author was supported by the
Centre for Research on Optimization and Control (CEOC) from the ``Funda\c{c}\~{a}o para a
Ci\^{e}ncia e a Tecnologia'' (FCT), cofinanced by the European Community Fund FEDER/POCTI.}}

\begin{document}

\maketitle

\begin{abstract}
We prove a multidimensional analogue of the classical Lagrange theorem on continued fractions. As a
multidimensional generalization of continued fractions we use Klein polyhedra.

%%Библиография: 15 наименований.
\end{abstract}

\section{Introduction}
\label{sec1}

The Lagrange theorem on continued fractions states that {\it a number~$\alpha$ is a quadratic surd
if and only if its continued fraction expansion is eventually periodic}.

We give a geometric interpretation of this fact. To this end, we consider two lines
in~$\mathbb{R}^2$, generated by vectors $(1,\alpha)$ and $(1,\beta)$,
$\alpha,\beta\in\mathbb{R}\setminus\mathbb{Q}$, $\alpha\neq\beta$. These lines divide the plane
into four angles, which we shall call {\it cones\/} minding the future multidimensional
generalizations. In the interior of each cone we consider the convex hull of integer points. The
obtained four unbounded convex polygons are called {\it Klein polygons}. There is a remarkable
correspondence between the partial quotients of the numbers~$\alpha$ and~$\beta$ and the edges of
the Klein polygons (see~\cite{1},~\cite{2}): the integer lengths of a Klein polygon's edges and the
integer angles between adjacent edges are equal to the respective partial quotients of the
numbers~$\alpha$ and~$\beta$. The {\it integer length\/} of a segment with endpoints
in~$\mathbb{Z}^2$ is defined as the number of lattice points contained in the interior of this
segment plus~1. And the {\it integer angle\/} between two such segments with a common endpoint is
defined as the area of the parallelogram spanned by them divided by the product of their integer
lengths. When the segments are non-parallel the integer angle is obviously equal to the index of
the sublattice spanned by the primitive integer vectors parallel to these two segments. The
described construction allows to give a geometric analogue of the Lagrange theorem (see~\cite{2}).

\begin{statement}
\label{st1} The following two statements are equivalent:

{\rm 1)} a cone~$\mathscr{C}$ is invariant under the action of some $\SL_2(\mathbb{Z})$--operator
with distinct real positive eigenvalues;

{\rm 2)} the combinatorial structure of the boundary of the Klein polygon related to a
cone~$\mathscr{C}$ equipped with integer lengths of edges and integer angles between adjacent
edges is periodic.
\end{statement}

A few years ago Vladimir Arnold proposed to find an analogue of this statement for so called {\it
Klein polyhedra}, a natural multidimensional generalization of Klein polygons.

\begin{defin}
\label{def1} Let $\mathscr{C}\subset\mathbb{R}^n$ be an $n$-dimensional simplicial cone with its
vertex in the origin~$\boldsymbol{0}$. The convex hull~$K$ of non-zero points of~$\mathbb{Z}^n$
contained in~$\mathscr{C}$ is called a {\it Klein polyhedron} related to the lattice~$\mathbb{Z}^n$
and the cone~$\mathscr{C}$.
\end{defin}

For generalization of Statement~\ref{st1} it is natural to consider the case when~$\mathscr{C}$ is
irrational, that is, when the planes spanned by the cone's faces do not contain lattice points
except~$\boldsymbol{0}$. Then, as was shown in~\cite{3}, the Klein polyhedron is a generalized
polyhedron, which means that its intersection with an arbitrary bounded polyhedron is itself a
polyhedron. Hence the boundary of~$K$ is in this case an $(n-1)$-dimensional polyhedral surface
homeomorphic to~$\mathbb{R}^{n-1}$, consisting of convex $(n-1)$-dimensional (generalized)
polyhedra. Some of these polyhedra may be unbounded, but anyway, each point of~$K$'s boundary
belongs only to finitely many faces of~$K$.

\begin{defin}
\label{def2} The boundary~$\Pi$ of a Klein polyhedron~$K$ is called a {\it sail}.
\end{defin}

Due to the correspondence described above we can consider faces of a sail and edge stars of its
vertices as multidimensional analogues of partial quotients. Vladimir Arnold conjectured that there
are local affine invariants of a sail, the periodicity of which implies the existence of a
nonidentity $\SL_n(\mathbb{Z})$-operator preserving the cone (and hence the sail). Soon after the
problem had been posed results started appearing in this direction (see~\cite{4}--\cite{8}).
Especially interesting is the paper~\cite{5}, in which the author proposes (without proof) a
multidimensional analogue of Statement~\ref{st1}. However, there is an inaccuracy in its
formulation. In the current paper we formulate accurately and prove the correspondent theorem
from~\cite{5}.

\section{Formulation of the main result}
\label{sec2}

For each vertex $\mathbf{v}$ of a sail~$\Pi$ we shall denote by~$\oSt_\mathbf{v}$ the edge star
of~$\mathbf{v}$, that is, the union of all the sail's edges incident to~$\mathbf{v}$.

An ordered set (finite or infinite) of vertices of a sail~$\Pi$ such that every two consequent
vertices are connected with an edge, will be called a {\it chain\/} of vertices of~$\Pi$. For each
positive integer~$k$ we consider the chains of length~$k$ of the sail's vertices and denote
by~$\mathscr{V}_k(\Pi)$ the set of all such chains. We consider also a graph~$\mathscr{G}_k(\Pi)$
with the set of vertices equal to the set~$\mathscr{V}_k(\Pi)$ and the set of edges equal to the
set of pairs $(V,W)\in\mathscr{V}_k(\Pi)\times\mathscr{V}_k(\Pi)$ such that $V\cup
W\in\mathscr{V}_{k+1}(\Pi)$. A graph $\mathscr{G}_{k+1}(\Pi)$ is obviously isomorphic to the edge
graph of~$\mathscr{G}_k(\Pi)$. As in the case of a sail's vertices, an ordered set (finite or
infinite) of vertices of~$\mathscr{G}_k(\Pi)$ such that every two consequent vertices are connected
with an edge, will be called a {\it chain\/} of vertices of~$\mathscr{G}_k(\Pi)$. It is obvious
that if a chain of vertices of~$\mathscr{G}_k(\Pi)$ has length~$l\geq 2$, then it corresponds
naturally to a chain of length~$l-1$ in~$\mathscr{G}_{k+1}(\Pi)$. In particular, each vertex chain
of length~$l$ in~$\Pi$ corresponds to a chain of length $l+1-k$ in~$\mathscr{G}_k(\Pi)$, for every
$k\leq l$.

Let us consider the group $\Aff_n(\mathbb{Z})$ of all the integer affine operators and a
subset~$\mathfrak{A}$ of this group. We define an {\it $\mathfrak{A}$-colouring\/}
of~$\mathscr{G}_k(\Pi)$ as follows. Two distinct vertices~$V$ and~$W$ of this graph are said to be
of the same colour if there is an operator in~$\mathfrak{A}$ that takes $\bigcup_{\mathbf{v}\in
V}\oSt_{\mathbf{v}}$ to~$\bigcup_{\mathbf{v}\in W}\oSt_{\mathbf{v}}$. We note that we do not
require this operator to preserve the order of the sail's vertices.

\begin{defin}
\label{def3} Let $\{\mathbf{v}_i\}_{i\in\mathbb{Z}}$ be a chain of vertices of a sail~$\Pi$.
Suppose that every $n-1$ consequent vertices in this chain are affinely independent (that is, do
not belong to any $(n-3)$-dimensional plane) and lie in an $(n-1)$-dimensional face of a sail.
Suppose that the images of this chain in~$\mathscr{G}_n(\Pi)$ and~$\mathscr{G}_{n+1}(\Pi)$ have
periodic $\mathfrak{A}$-colourings, and suppose also that for any two distinct operators
$A,B\in\mathfrak{A}$ ``establishing'' the colouring, the operator $AB^{-1}$ also belongs
to~$\mathfrak{A}$. Then we say that $\{\mathbf{v}_i\}_{i\in\mathbb{Z}}$ is {\it
$\mathfrak{A}$-periodic}.
\end{defin}

Let us denote by~$\mathfrak{A}_0$ the set of affine operators~$\mathbf{A}$ such that
$\mathbf{A}\colon\mathbf{x}\mapsto A(\mathbf{x})+\mathbf{a}$, $A\in\SL_n(\mathbb{Z})$,
$\mathbf{a}\in\mathbb{Z}^n$, and~$A$ satisfies the following two conditions:

(P1) all the eigenvalues of~$A$ are different from the unit;

(P2) if $\alpha\in\mathbb{C}\setminus\mathbb{R}$ is an eigenvalue of~$A$, then all the other
eigenvalues except for the complex conjugate of~$\alpha$ have absolute values different from that
of~$\alpha$.
\smallskip

The following theorems (especially Theorems~\ref{th1},~\ref{th4} and~\ref{th5}) refine the result
of the paper~\cite{5}.

\begin{Theorem}
\label{th1} Given an irrational cone $\mathscr{C}\subset\mathbb{R}^n$ consider the sail~$\Pi$
corresponding to~$\mathscr{C}$ and~$\mathbb{Z}^n$. The following two statements are equivalent:

{\rm 1)} there is a nonidentity operator $A\in\SL_n(\mathbb{Z})$ such that
$A(\mathscr{C})=\mathscr{C}$;

{\rm 2)} there is an unbounded in both directions (as a subset of~$\mathbb{R}^n$)
$\mathfrak{A}_0$-periodic chain $\{\mathbf{v}_i\}_{i\in\mathbb{Z}}$ of vertices of~$\Pi$.
\end{Theorem}

The implication 1)$\Rightarrow$2) is obvious since if $A(\mathscr{C})=\mathscr{C}$, then
$A(\Pi)=\Pi$. The implication 2)$\Rightarrow$1) follows from Theorems~\ref{th2} and~\ref{th3},
which are proved in the current paper.

\begin{Theorem}
\label{th2} Let $\mathfrak{A}$ be an arbitrary subset of $\Aff_n(\mathbb{Z})$ and let
$\{\mathbf{v}_i\}_{i\in\mathbb{Z}}$ be an $\mathfrak{A}$-periodic chain of vertices of a sail
$\Pi\subset\mathbb{R}^n$. Then there is an operator in~$\mathfrak{A}$ establishing a nontrivial
shift of the set $\bigcup_{i\in\mathbb{Z}}\oSt_{\mathbf{v}_i}$ along itself.
\end{Theorem}

\begin{defin}
\label{def4} An operator $A\in\SL_n(\mathbb{Z})$ is called {\it hyperbolic} if its characteristic
polynomial is irreducible over~$\mathbb{Q}$ with all the roots real and positive.
\end{defin}

\begin{Theorem}
\label{th3} Let a sail $\Pi$ correspond to the lattice~$\mathbb{Z}^n$ and an irrational cone
$\mathscr{C}\subset\mathbb{R}^n$. Let $\{\mathbf{v}_i\}_{i\in\mathbb{Z}}$ be an unbounded in both
directions (as a subset of~$\mathbb{R}^n$) sequence (not necessarily a chain) of vertices of~$\Pi$.
Suppose that there is an operator~$\mathbf{A}\in\mathfrak{A}_0$ such that
$\mathbf{A}(\oSt_{\mathbf{v}_i})=\oSt_{\mathbf{v}_{i+1}}$. Then $\mathbf{A}=A\in\SL_n(\mathbb{Z})$
and $A(\mathscr{C})=\mathscr{C}$. Besides that, if all the eigenvalues of~$A$ are pairwise
distinct, then~$A$ is a hyperbolic operator.
\end{Theorem}

\begin{Remark}
\label{rem1} It follows from the proof of Theorem~\ref{th3} that if the operator
    establishing the shift of the chain is linear, then we may confine our requirements to
    the property~(P2) and neglect the property~(P1).
\end{Remark}

Let us denote by~$\mathfrak{A}_1$ the set of operators from~$\mathfrak{A}_0$ with linear component
having pairwise distinct eigenvalues. Then we obtain another theorem, which is also a corollary of
Theorems~\ref{th2} and~\ref{th3}.

\begin{Theorem}
\label{th4} Let a sail~$\Pi$ correspond to the lattice~$\mathbb{Z}^n$ and an irrational cone
$\mathscr{C}\subset\mathbb{R}^n$. Then the following two statements are equivalent:

{\rm 1)} there is a hyperbolic operator $A\in\SL_n(\mathbb{Z})$ such that
$A(\mathscr{C})=\mathscr{C}$;

{\rm 2)} there is an unbounded in both directions (as a subset of~$\mathbb{R}^n$)
$\mathfrak{A}_1$-periodic chain $\{\mathbf{v}_i\}_{i\in\mathbb{Z}}$ of vertices of~$\Pi$.
\end{Theorem}

It is worth mentioning that due to the hyperbolicity of~$A$ the implication 1)$\Rightarrow$2)
in~Theorem~\ref{th4} also follows from the Dirichlet theorem on algebraic units (see~\cite{9}
and~\cite{10}). We also note that in~Theorem~\ref{th1} one cannot keep~$\mathfrak{A}_0$ and add the
requirement of hyperbolicity to the statement~1). For even~$n$ it is easy to find an operator
$A\in\SL_n(\mathbb{Z})$ having an irrational invariant (simplicial) cone and the characteristic
polynomial equal, for instance, to the square of a polynomial with integer coefficients and
positive real roots, irreducible over~$\mathbb{Q}$. For~$n=\nobreak4$ one can consider the operator
with the matrix
$$
A=\begin{pmatrix}
2 & 1 & 0 & 0
\\
1 & 1 & 0 & 0
\\
0 & 0 & 2 & 1
\\
0 & 0 & 1 & 1
\end{pmatrix}.
$$
This operator has two invariant irrational two-dimensional planes; in each of these planes we can
choose two vectors in such a way that the four vectors thus obtained generate an irrational
cone~$\mathscr{C}$, invariant under the action of~$A$, but at the same time not invariant under the
action of any hyperbolic integer operator.

On the other hand one cannot add the requirement of hyperbolicity in the statement~1) of
Theorem~\ref{th1} having only replaced the property~(P2) with the property of an operator to have
pairwise distinct eigenvalues. It can be seen from the following example kindly provided by Elena
Korkina (Pavlovskaia).

\begin{Example}
Consider a quadratic equation $\lambda^2-p\lambda+1=0$ with integer $p\geq 3$. It has two distinct
real irrational positive roots~$\lambda_1$ and~$\lambda_2$. After replacing~$\lambda$ with~$\mu^3$
we get an equation of $6$-th degree: $\mu^6-p\mu^3+1=0$. Consider an extension of~$\mathbb{Q}$ with
a root of this equation. It is a $6$-dimensional vector space~$M$ over~$\mathbb{Q}$ with the
basis~$1,\mu,\dots,\mu^5$. Consider an $\SL_6(\mathbb{Z})$-operator~$A$ acting on~$M$ as
multiplication by~$\mu$. Let us imbed naturally~$M$ into~$\mathbb{R}^6$ and extend the action
of~$A$ from the image of~$M$ to~$\mathbb{R}^6$.

The operator~$A$ has six eigenvalues: $\sqrt[3]{\lambda_1}$\,, $\zeta\sqrt[3]{\lambda_1}$\,,
$\zeta^2\sqrt[3]{\lambda_1}$ and $\sqrt[3]{\lambda_2}$\,, $\zeta\sqrt[3]{\lambda_2}$\,,
$\zeta^2\sqrt[3]{\lambda_2}$\,, where $\zeta=e^{2\pi i/3}$. For each $i=1,2$ let us denote by~$V_i$
the three-dimensional $A$-invariant subspace of~$\mathbb{R}^6$ corresponding to the eigenvalues
$\sqrt[3]{\lambda_i}$\,, $\zeta\sqrt[3]{\lambda_i}$\,, $\zeta^2\sqrt[3]{\lambda_i}$ and consider an
arbitrary vector $\boldsymbol{\omega}_i\in V_i$ not belonging to any one- or two-dimensional
invariant subspace of~$A$. Let $\boldsymbol{\omega}_i^{0}=\boldsymbol{\omega}_i$,
$\boldsymbol{\omega}_i^{1}=A\boldsymbol{\omega}_i$,
$\boldsymbol{\omega}_i^{2}=A^2\boldsymbol{\omega}_i$. The vectors $\boldsymbol{\omega}_1^{0}$,
$\boldsymbol{\omega}_1^{1}$, $\boldsymbol{\omega}_1^{2}$, $\boldsymbol{\omega}_2^{0}$,
$\boldsymbol{\omega}_2^{1}$, $\boldsymbol{\omega}_2^{2}$ are obviously linearly independent.
Consider the cone~$\mathscr{C}$ generated by these vectors. Its three-dimensional faces generated
by the first three and the last three vectors are invariant under the action of~$A$, since
$A\boldsymbol{\omega}_i^{2}=A^3\boldsymbol{\omega}_i^{0}= \lambda_i\boldsymbol{\omega}_i^{0}$ and
$\lambda_i>0$. Hence the cone~$\mathscr{C}$ and the sail corresponding to it are invariant under
the action of~$A$. However, the operator $A$ is not hyperbolic. Moreover, one can choose the
vectors~$\boldsymbol{\omega}_1$ and~$\boldsymbol{\omega}_2$ to generate transcendental directions,
thus preventing~$\mathscr{C}$ from being invariant under the action of any hyperbolic integer
operator.
\end{Example}

For $n=3$ a stronger statement than Theorems~\ref{th1} and~\ref{th4} can be proved. We can
replace~$\mathfrak{A}_0$ by $\Aff_3(\mathbb{Z})$ keeping the requirement of hyperbolicity. It is
also clear that if $n=3$ and $\mathfrak{A}=\Aff_3(\mathbb{Z})$ the technical conditions of
Definition~\ref{def3} are automatically satisfied and $\Aff_3(\mathbb{Z})$-periodicity of a chain
of vertices means exactly the periodicity of the $\Aff_3(\mathbb{Z})$-colouring of the chain's
images in~$\mathscr{G}_n(\Pi)$ and~$\mathscr{G}_{n+1}(\Pi)$.

\begin{Theorem}
\label{th5} Let a sail~$\Pi$ correspond to the lattice~$\mathbb{Z}^3$ and an irrational cone
$\mathscr{C}\subset\mathbb{R}^3$. Then the following two statements are equivalent:

{\rm 1)} there is a hyperbolic operator $A\in\SL_3(\mathbb{Z})$ such that
$A(\mathscr{C})=\mathscr{C}$;

{\rm 2)} there is an unbounded in both directions (as a subset of~$\mathbb{R}^3$)
$\Aff_3(\mathbb{Z})$-periodic chain $\{\mathbf{v}_i\}_{i\in\mathbb{Z}}$ of vertices of~$\Pi$.
\end{Theorem}

\section{Proof of Theorem~\ref{th2}}
\label{sec3}

Let $\{V_i\}_{i\in\mathbb{Z}}$ and $\{W_i\}_{i\in\mathbb{Z}}$ denote the images of
$\{\mathbf{v}_i\}_{i\in\mathbb{Z}}$ in~$\mathscr{G}_n(\Pi)$ and $\mathscr{G}_{n+1}(\Pi)$,
respectively. Suppose that $V_0\setminus V_1=W_0\setminus W_1=\{\mathbf{v}_0\}$. Set
$$
\widetilde V_i=\bigcup_{\mathbf{v}\in V_i}\oSt_{\mathbf{v}},\qquad
\widetilde W_i=\bigcup_{\mathbf{v}\in W_i}\oSt_{\mathbf{v}}.
$$
Consider the periods of the colourings of $\mathscr{G}_n(\Pi)$ and~$\mathscr{G}_{n+1}(\Pi)$.
Set~$t$ to be equal to their least common multiple. Then there are operators
$\mathbf{A},\mathbf{B},\mathbf{C}\in\mathfrak{A}$ such that
$$
\mathbf{A}(\widetilde V_0)=\widetilde V_t,\quad
\mathbf{B}(\widetilde V_1)=\widetilde V_{t+1},\quad
\mathbf{C}(\widetilde W_0)=\widetilde W_t.
$$

Suppose that~$\mathbf{A}$ inverts the order of the vertices, that is,
$\mathbf{A}(\mathbf{v}_i)\,{=}\,\mathbf{v}_{t+n-1-i}$ for all $i=0,1,\dots,n-1$. Then, either the
operator~$\mathbf{A}'$ such that $\mathbf{A}'(\widetilde V_t)=\widetilde V_{2t}$, or the operator
$\mathbf{A}'\mathbf{A}$ (for which, obviously, $\mathbf{A}'\mathbf{A}(\widetilde V_0)=\widetilde
V_{2t}$) preserves the order of the vertices. Therefore we may assume that~$\mathbf{A}$ preserves
the order of the vertices, doubling~$t$, if necessary, or shifting all the indices by~$t$.

Suppose now that~$\mathbf{B}$ inverts the order of the vertices. Then it is easy to see that
either~$\mathbf{C}$, or $\mathbf{B}\mathbf{C}^{-1}\mathbf{A}\mathbf{C}^{-1}\mathbf{B}$
takes~$\widetilde V_1$ to~$\widetilde V_{t+1}$ preserving the order of the vertices. Therefore we
may assume that~$\mathbf{B}$ also preserves the order of the vertices.

Assuming that~$\mathbf{A}$ and~$\mathbf{B}$ preserve the order of the vertices, let us consider the
operator $\mathbf{B}^{-1}\mathbf{A}$. It is obvious that
$\mathbf{B}^{-1}\mathbf{A}(\mathbf{v}_i)=\mathbf{v}_i$ and
$\mathbf{B}^{-1}\mathbf{A}(\oSt_{\mathbf{v}_i})=\oSt_{\mathbf{v}_i}$ for each $i=1,2,\dots,n-1$.
Let~$\mathbf{r}_0$ denote the sum of primitive vectors parallel to the edges
of~$\oSt_{\mathbf{v}_1}$, and let $\mathbf{r}_i=\mathbf{v}_{i+1}-\mathbf{v}_i$ for
$i=1,2,\dots,n-2$. It is clear that $\mathbf{B}^{-1}\mathbf{A}(\mathbf{v}_1+\mathbf{r}_i)=
\mathbf{v}_1+\mathbf{r}_i$ for every $i=0,1,\dots,n-2$. It follows from the assumption of the
theorem that $\mathbf{r}_1,\dots,\mathbf{r}_{n-2}$ are linearly independent and parallel to some
$(n-1)$-dimensional face, incident to~$\mathbf{v}_1$. Hence the vectors
$\mathbf{r}_0,\mathbf{r}_1,\dots,\mathbf{r}_{n-2}$ are also linearly independent. Thus, the affine
hull of the points $\mathbf{v}_1,\mathbf{v}_1+\mathbf{r}_0,\mathbf{v}_1+\mathbf{r}_1,\dots,
\mathbf{v}_1+\mathbf{r}_{n-2}$ is an $(n-1)$-dimensional invariant plane of the operator
$\mathbf{B}^{-1}\mathbf{A}$. This plane divides the edges of~$\oSt_{\mathbf{v}_1}$ into two
invariant subsets, since $\mathbf{B}^{-1}\mathbf{A}$ preserves the orientation. Let~$\mathbf{r}_0'$
and~$\mathbf{r}_0''$ denote the sums of primitive vectors parallel to the edges from these two
subsets. Then $\mathbf{B}^{-1}\mathbf{A}$ preserves $\mathbf{v}_1+\mathbf{r}_0'$ and
$\mathbf{v}_1+\mathbf{r}_0''$ and at least one of these points does not belong to the affine hull
of $\mathbf{v}_1,\mathbf{v}_1+\mathbf{r}_0,\mathbf{v}_1+\mathbf{r}_1,\dots,
\mathbf{v}_1+\mathbf{r}_{n-2}$. We get $n+1$ points, invariant under the action of
$\mathbf{B}^{-1}\mathbf{A}$, which do not belong simultaneously to any $n$-dimensional plane. Hence
$\mathbf{B}^{-1}\mathbf{A}$ is an identity operator and $\mathbf{B}=\mathbf{A}$.

Continuing these arguments in both directions of the chain $\{\mathbf{v}_i\}_{i\in\mathbb{Z}}$ we
get that~$\mathbf{A}$ shifts $\bigcup_{i\in\mathbb{Z}}\oSt_{\mathbf{v}_i}$ along itself.

\section{A theorem on an edge star and integer distance}
\label{sec4}

In this section we formulate and prove Theorem~\ref{th6}, which claims that if one can ``see'' a
convex integer polyhedron from a point of~$\mathbb{Z}^n$, then this point cannot be too ``far''
from the polyhedron, and the ``distance'' is bounded from above by a constant depending only on the
integer-combinatorial structure of the edge stars of the polyhedron. We shall use this theorem to
prove Theorem~\ref{th3} and also in the end of this paper to reformulate the Oppenheim conjecture.

Let $\conv(M)$ and $\aff(M)$ denote respectively the convex and the affine hulls of a set
$M\subset\mathbb{R}^n$.

\begin{defin}
\label{def5} A (generalized) polyhedron with all its vertices in~$\mathbb{Z}^n$ is called {\it
integer}.
\end{defin}

\begin{defin}
\label{def6} Suppose that a vertex~$\mathbf{v}$ of an integer convex $n$-dimensional (generalized)
polyhedron~$P$ is incident to~$m$ edges. Let $\mathbf{r}_1,\dots,\mathbf{r}_m$ denote the primitive
vectors of~$\mathbb{Z}^n$ parallel to these edges. Then we define the {\it determinant\/} of the
edge star~$\oSt_{\mathbf{v}}$ as
$$
\det\oSt_{\mathbf{v}}=\sum_{1\leq i_1<\dots<i_n\leq m}
|\det(\mathbf{v}_{i_1},\dots,\mathbf{v}_{i_n})|.
$$
In other words, $\det\oSt_{\mathbf{v}}$ is equal to the volume of the Minkowski sum of the segments
$[\boldsymbol{0},\mathbf{r}_1],\dots,[\boldsymbol{0},\mathbf{r}_m]$.
\end{defin}

\begin{defin}
\label{def7} Let $F\subset\mathbb{R}^n$ be an $(n-1)$-dimensional integer polyhedron and
$\mathbf{a}\in\mathbb{Z}^n$. The index of the minimal (with respect to inclusion) sublattice
of~$\mathbb{Z}^n$ containing the set $\mathbb{Z}^n\cap\aff(F)-\mathbf{a}$ is called an {\it integer
distance\/} from~$F$ to~$\mathbf{a}$ and is denoted by $\rho_{\opint}(F,\mathbf{a})$.
\end{defin}

\begin{Theorem}
\label{th6} Let $P\subset\mathbb{R}^n$ be an $n$-dimensional integer convex polyhedron, not
containing~$\boldsymbol{0}$. Let $\mathbf{v}$ be a vertex of~$P$ contained in the interior of
$\conv(P\cup\{\boldsymbol{0}\})$, let~$\mathbf{v}$ be incident to~$m$ edges of~$P$, and let $F$ be
an $(n-1)$-dimensional face of~$P$ incident to~$\mathbf{v}$. Let also
$$
\mathbb{Z}^n\cap\conv(P\cup\{\boldsymbol{0}\})\setminus P=\{\boldsymbol{0}\}.
$$
Then
$$
\rho_{\opint}(F,\boldsymbol{0})<(nm)^{4n!}\det\oSt_{\mathbf{v}}.
$$
\end{Theorem}

When proving Theorem~\ref{th3} we shall actually use a corollary of Theorem~\ref{th6} rather than
Theorem~\ref{th6} itself.

\begin{corollary}
\label{cor1} Let $\mathbf{v}$ be a vertex of a sail~$\Pi$ and let $F$ be an $(n-1)$-dimensional
face of~$\Pi$ incident to~$\mathbf{v}$. Then
$$
\rho_{\opint}(F,\boldsymbol{0})<(n\det\oSt_{\mathbf{v}})^{8n!}.
$$
\end{corollary}

We also mention another corollary of Theorem~\ref{th6}, which can be useful for an explicit
description of a Klein polyhedron's faces, though we do not make any use of it in the current
paper.

\begin{corollary}
\label{cor2} Let $P\subset\Delta\subset\mathbb{R}^n$ be an $n$-dimensional integer convex
polyhedron contained in an arbitrary polyhedron~$\Delta$ (for instance, in a simplex), and let
$$
P=\bigcap_{i=1}^k\{\mathbf{x}\in\mathbb{R}^n \mid
\langle\mathbf{x},\mathbf{f}_i\rangle\geq D_i\},
$$
where $\mathbf{f}_i$ are primitive vectors of the lattice~$\mathbb{Z}^n$. Let every vertex of~$P$
be incident to not more than~$m$ edges, let $D$ be a constant bounding from above the determinants
of all the edge stars of~$P$ and let
$$
\mathbb{Z}^n\cap\bigcap_{i=1}^k\bigl\{\mathbf{x}\in\mathbb{R}^n\mid
\langle\mathbf{x},\mathbf{f}_i\rangle>D_i-D(nm)^{4n!}\bigr\}=
\mathbb{Z}^n\cap P.
$$
Then $\mathbb{Z}^n\cap\Delta=\mathbb{Z}^n\cap P$.
\end{corollary}

To prove Theorem~\ref{th6} we shall need a few auxiliary statements.

\begin{Lemma}
\label{lem1} Let $n$ be a positive integer and let~$V$,~$A$, $B$ be real numbers such that $V>0$,
$A>V/n$, $0<B<1$. Let $f(x)=A(B-x)^{n-1}(1-x)-x$. Then $f(x)>0$ for all~$x$ satisfying the
inequality $0\leq x\leq B^{n-1}n^{-1}(1+V^{-1})^{-1}$.
\end{Lemma}

\begin{proof}
The root of the equation of the tangent line for~$f(x)$ in zero equals
$$
\frac{AB^{n-1}}{1+AB^{n-2}(n-1+B)}>\frac{B^{n-1}}{n(1+V^{-1})}
$$
and bounds from below the minimal positive root of~$f(x)$.
\end{proof}

\begin{Lemma}
\label{lem2} Let $\Lambda$ be an $n$-dimensional lattice in~$\mathbb{R}^n$ with determinant equal
to~$1$. Let $\Delta\subset\mathbb{R}^n$ be an $n$-dimensional simplex with vertices
$\mathbf{v}_0\in\mathbb{R}^n$ and $\mathbf{v}_1,\dots,\mathbf{v}_n\in\mathbb{Q}\Lambda$. Let a
point $\mathbf{v}\in\Lambda$ belong to its interior. Let~$\Delta_0$ denote the simplex with
vertices $\mathbf{v},\mathbf{v}_1,\dots,\mathbf{v}_n$. Suppose that there is a constant~$V$ such
that for every $k\in\{1,\dots,n\}$ the $k$-dimensional volume of every $k$-dimensional face
of~$\Delta_0$ is at least $V$ times greater than the determinant of the $k$-dimensional (affine)
lattice consisting of the points of~$\Lambda$ contained in the affine hull of this face. Suppose
also that
$$
\mathbb{Z}^n\cap\biggl(\mathbf{v}_0+\biggl(1-\frac{\opvol_n(\Delta_0)}
{\opvol_n(\Delta)}\biggr)(\Delta-\mathbf{v}_0)\biggr)\Big\backslash
\{\mathbf{v},\mathbf{v}_0\}=\varnothing
$$
(that is, the hyperplane containing~$\mathbf{v}$ and parallel to the plane
$\aff(\mathbf{v}_1,\dots,\mathbf{v}_n)$ cuts off a simplex from~$\Delta$, which does not contain
any lattice points, except for~$\mathbf{v}$ and maybe $\mathbf{v}_0$). Then
$$
\opvol_n(\Delta)\leq\opvol_n(\Delta_0)\prod_{k=1}^nk^\frac{n!}{k!}
(1+V^{-1})^\frac{(n-1)!}{(k-1)!}.
$$
\end{Lemma}

\begin{proof}
Let us set
$$
c_n(V)=\prod_{k=1}^nk^\frac{n!}{k!}(1+V^{-1})^\frac{(n-1)!}{(k-1)!}
$$
and apply induction by~$n$. For $n=1$ the statement of the lemma is obvious. Suppose that the
statement of the lemma is true in dimension $n-1$. Let us prove it for dimension~$n$. Due to the
affinity of the problem, it is sufficient to prove the statement of the lemma in the case when the
vectors $\mathbf{v}_1-\mathbf{v}_0,\dots,\mathbf{v}_n-\mathbf{v}_0$ are pairwise orthogonal and
have absolute values equal to some number~$l$.

Without loss of generality we assume that among the points $\mathbf{v}_1,\dots,\mathbf{v}_n$ the
point $\mathbf{v}_n$ is farthest from~$\mathbf{v}$. This means that among the faces
$\conv(\{\mathbf{v}_0,\mathbf{v}_1,\dots,\mathbf{v}_n\}\setminus \{\mathbf{v}_i\})$ the face
$\conv(\mathbf{v}_0,\mathbf{v}_1,\dots,\mathbf{v}_{n-1})$ is closest to~$\mathbf{v}$. Suppose that
$\mathbf{u}\in\conv(\mathbf{v}_0,\mathbf{v}_1,\dots,\mathbf{v}_{n-1})$ is the closest point
to~$\mathbf{v}$ such that $\mathbf{u}-\mathbf{v}$ is parallel to the plane
$\pi=\aff(\mathbf{v}_1,\dots,\mathbf{v}_n)$. Then the set
$$
\Omega=\conv(\{\mathbf{v}_0\}\cup
(\pi\cap B_{|\mathbf{u}-\mathbf{v}|}(\mathbf{v}))),
$$
where $B_{|\mathbf{u}-\mathbf{v}|}(\mathbf{v})$ is the ball of radius $|\mathbf{u}-\mathbf{v}|$
centered at~$\mathbf{v}$, does not contain any points of~$\Lambda$, except for~$\mathbf{v}$ and
maybe $\mathbf{v}_0$. Consequently, the volume of~$\Omega$, due to the Minkowski theorem on convex
bodies, does not exceed $2^{n-1}$. This means that
\begin{equation}
\label{eq1}
\frac{|\mathbf{u}-\mathbf{v}|^{n-1}}{n}
\biggl(1-\frac{\opvol_n(\Delta_0)}{\opvol_n(\Delta)}\biggr)
\frac{l}{\sqrt n}\leq2^{n-1}\frac{\opvol_n(\Delta_0)}{V}\,.
\end{equation}

Set
\begin{align*}
\Delta_0'&=\conv(\mathbf{v},\mathbf{v}_1,\dots,\mathbf{v}_{n-1}),
\\
\Delta'&=\aff(\Delta_0')\cap\Delta.
\end{align*}
The vertices of the simplex~$\Delta'$ are the points $\mathbf{v}_1,\dots,\mathbf{v}_{n-1}$ and some
point~$\mathbf{v}_0'$ from the segment $[\mathbf{v}_0,\mathbf{v}_n]$. By the induction assumption,
$\opvol_{n-1}(\Delta')\leq c_{n-1}(V)\*\opvol_{n-1}(\Delta_0')$, which, due to elementary geometric
consideration, implies that
\begin{equation}
\label{eq2}
|\mathbf{u}-\mathbf{v}|\geq\sqrt n\biggl(c_{n-1}^{-1}(V)\frac{l}{\sqrt{n-1}}-
\frac{\opvol_n(\Delta_0)}{\opvol_n(\Delta)}\,\frac{l}{\sqrt{n-1}}\biggr).
\end{equation}

Using~\eqref{eq1},~\eqref{eq2} and the fact that $\opvol_n(\Delta)=l^n/n!$ we get the inequality
$$
\frac{V(n-1)!(\sqrt{n-1})^{n-1}}{2^{n-1}\sqrt n(\sqrt{n-1})^{n-1}}
\biggl(c_{n-1}^{-1}(V)
-\frac{\opvol_n(\Delta_0)}{\opvol_n(\Delta)}\biggr)^{n-1}
\biggl(1-\frac{\opvol_n(\Delta_0)}{\opvol_n(\Delta)}\biggr)-
\frac{\opvol_n(\Delta_0)}{\opvol_n(\Delta)}\leq 0.
$$
Hence, due to Lemma~\ref{lem1},
$$
\frac{\opvol_n(\Delta_0)}{\opvol_n(\Delta)}>
\frac{c_{n-1}^{1-n}(V)}{n(1+V^{-1})}=c_n^{-1}(V),
$$
which completes the proof.
\end{proof}

We shall denote by $\opint P$ and $\opext P$ the relative interior and the vertex set of a
polyhedron~$P$. If $M\subset\mathbb{R}^n$ is a finite set and to each point $\mathbf{x}\in M$ a
positive mass~$\nu_{\mathbf{x}}$ is assigned, then for each subset $M'\subseteq M$ of
cardinality~$\sharp(M')$ we shall denote by~$c(M')$ the point $(\sum_{\mathbf{x}\in
M'}\nu_{\mathbf{x}}\mathbf{x})/\sharp(M')$, that is, the center of mass of the set~$M'$.

\begin{Lemma}
\label{lem3} Let $P$ be a convex $(n-1)$-dimensional polyhedron with arbitrary positive masses
assigned to its vertices. Let~$\mathfrak T$ be an arbitrary partition of the (relative) boundary
of~$P$ into (closed) simplices with vertices in $\opext P$. Then
$$
\opint P=\bigcup_{\Delta\in\mathfrak T}
\opint(\conv(\Delta\cup\{c(\opext P\setminus\opext\Delta)\})).
$$
\end{Lemma}

\begin{proof}
Let $\mathbf{x}\in\opint P$. It is obvious that there is a simplex $\Delta\in\mathfrak T$ such that
$\mathbf{x}\in\conv(\Delta\cup\{c(P)\})$. It remains to notice that
$\conv(\Delta\cup\{c(P)\})\cap\opint P$ is contained in the interior of $\conv(\Delta\cup\{c(\opext
P\setminus\opext\Delta)\})$.
\end{proof}

\begin{nproof}[Proof of Theorem~\ref{th6}]
Let $\mathbf{r}_1,\dots,\mathbf{r}_m$ be the primitive vectors of~$\mathbb{Z}^n$ parallel to the
edges incident to~$\mathbf{v}$. Take arbitrary positive numbers $k_1,\dots,k_m$ such that the
points $\mathbf{r}'_i=k_i\mathbf{r}_i$ lie on a same hyperplane, and set
$P=\conv(\mathbf{r}'_1,\dots,\mathbf{r}'_m)$. Consider a number~$\lambda$ such that
$\lambda\mathbf{v}\in P$.

Assign masses~$k_i^{-1}$ to the points $r'_i$. Then, by Lemma~\ref{lem3} we can renumerate the
vectors $\mathbf{r}_1,\dots,\mathbf{r}_m$ (changing respectively numerations of numbers
$k_1,\dots,k_m$ and vectors $\mathbf{r}'_1,\dots,\mathbf{r}'_m$) so that
$\lambda\mathbf{v}=\lambda'_0\mathbf{r}'_0+\nobreak\dots+\nobreak \lambda'_{n-1}\mathbf{r}'_{n-1}$
with strictly positive~$\lambda'_i$ and $\mathbf{r}'_0=(\mathbf{r}_n+\dots+\mathbf{r}_m)/(m-n+1)$.
Setting $\mathbf{r}_0=\mathbf{r}'_0,\ \lambda_0=\lambda'_0$ and $\lambda_i=k_i\lambda'_i$ for
$i=1,\dots,n-1$ we get that
$\lambda\mathbf{v}=\lambda_0\mathbf{r}_0+\dots+\lambda_{n-1}\mathbf{r}_{n-1}$ with strictly
positive~$\lambda_i$. The point $\mathbf{v}$ is therefore contained in the interior of the
simplex~$\Delta$ with vertices
$\boldsymbol{0},\mathbf{v}+\mathbf{r}_0,\dots,\mathbf{v}+\mathbf{r}_{n-1}$. It follows from the
definition of~$\mathbf{r}_0$ that for every $k\in\{1,\dots,n\}$ the $k$-dimensional volume of every
$k$-dimensional face of the simplex $\Delta_0=\conv(\mathbf{v},\mathbf{v}+\mathbf{r}_0,\dots,
\mathbf{v}+\mathbf{r}_{n-1})$ is at least $n^{-1}(m-n+1)^{-1}$ times greater than the determinant
of the $k$-dimensional lattice of integer points contained in the affine hull of this face.
Moreover, it is clear that the hyperplane containing~$\mathbf{v}$ and parallel to
$\aff(\mathbf{v}+\mathbf{r}_0,\dots,\mathbf{v}+\mathbf{r}_{n-1})$ is a support plane for the
polyhedron~$P$ and therefore cuts off a simplex from~$\Delta$ with no integer points inside it,
other than~$\boldsymbol{0}$ and~$\mathbf{v}$. Hence~$\Delta$ and~$\Delta_0$ satisfy the conditions
of Lemma~\ref{lem2}, which means that
\begin{align*}
\opvol_n(\Delta)&\leq\opvol_n(\Delta_0)\prod_{k=1}^nk^\frac{n!}{k!}
(mn-n^2+n+1)^\frac{(n-1)!}{(k-1)!}
\\
&<\opvol_n(\Delta_0)\prod_{k=1}^n(nmk)^\frac{n!}{k!}<
(nm)^{4n!}\opvol_n(\Delta_0).
\end{align*}
Therefore,
$$
\rho_{\opint}(F,\boldsymbol{0})\leq n\opvol_n(\conv(F\cup\{\boldsymbol{0}\}))
<n\opvol_n(\Delta)<(nm)^{4n!}\det\oSt_{\mathbf{v}},
$$
which proves the theorem.
\end{nproof}

\section{Proof of Theorem~\ref{th3}}
\label{sec5}

To prove Theorem~\ref{th3} we need the following fact, which we give without proof because of its
simplicity.

\begin{Lemma}
\label{lem4} Let $A\in\GL_n(\mathbb{R})$ be an $n$-dimensional Jordan cell with its eigenvalue
equal to~$\lambda$. Let $\mathbf{e}\in\mathbb{R}^n$ be a vector with the $i$-th coordinate equal
to~$1$ and all the other coordinates equal to zero. Let $\mathbf{h}\in\mathbb{R}^n$ be a vector
with nonzero coordinates. Then $\langle A^m\mathbf{h},\mathbf{e}\rangle\asymp\lambda^mm^{n-i}$ as
$m\to+\infty$.
\end{Lemma}

\begin{nproof}[Proof of Theorem~\ref{th3}]
Due to~(P1) the operator $E-A$ is invertible. Set $\mathbf{b}=(E-A)^{-1}\mathbf{a}$. Then
$\mathbf{b}\in\mathbb{Z}^n$ and $\mathbf{A}(\mathbf{b}+\mathbf{x})=\mathbf{b}+A(\mathbf{x})$. Thus,
we can assume that instead of an affine operator~$\mathbf{A}$ we have a linear operator~$A$ acting,
but the vertex of~$\mathscr{C}$ is in~$-\mathbf{b}$.

If $F$ is an arbitrary $(n-1)$-dimensional face of~$\Pi$ incident to~$\mathbf{v}_0$, then by
Corollary~\ref{cor1} the integer distance from~$A^m(F)$ to~$-\mathbf{b}$ is bounded for all
$m\in\mathbb{Z}$ by a constant not depending on~$m$. Hence, if
$\langle\mathbf{h}_F,\,\cdot\,\rangle$ is a linear form such that
$\aff(F)=\{\mathbf{x}\in\mathbb{R}^n\mid \langle\mathbf{h}_F,\mathbf{x}\rangle=1\}$, then
\begin{equation}
\label{eq3}
\max_{m\in\mathbb{Z}}|\langle(A^\ast)^m\mathbf{h}_F,-\mathbf{b}\rangle|
<\infty.
\end{equation}

Consider a Jordan basis $\mathscr{E}=\{\mathbf{e}_1,\dots,\mathbf{e}_n\}$ of~$A$. Consider also a
linear form $\langle\mathbf{h},\,\cdot\,\rangle$ such that neither of the components
of~$\mathbf{h}$ in the basis~$\mathscr{E}^\ast$ dual to~$\mathscr{E}$ is equal to zero and such
that the plane $\pi_{\mathbf{h}}=\{\mathbf{x}\in\mathbb{R}^n\mid
\langle\mathbf{h},\mathbf{x}\rangle=1\}$ intersects the edges of~$\oSt_{\mathbf{v}_0}$ by inner
points. It follows from~\eqref{eq3} that
$$
\max_{m\in\mathbb{Z}}|\langle(A^\ast)^m\mathbf{h},-\mathbf{b}\rangle|<\infty.
$$
This, together with Lemma~\ref{lem4} applied in the cases $m\to+\infty$ and $m\to-\infty$, implies
that
\begin{equation}
\label{eq4}
|A(\mathbf{b})|=|\mathbf{b}|.
\end{equation}

On the other hand, the fact that the plane~$\pi_{\mathbf{h}}$ cuts the edges
of~$\oSt_{\mathbf{v}_0}$ by inner points implies that
\begin{equation}
\label{eq5}
\inf_{m\in\mathbb{Z}}\opvol_n(\conv(-\mathbf{b},\mathscr{C}\cap
A^m(\pi_{\mathbf{h}})))>0.
\end{equation}

Let $\boldsymbol{\omega}_1,\dots,\boldsymbol{\omega}_n$ denote the unit vectors generating the
edges of~$\mathscr{C}$. Then
\begin{align}
\nonumber
A^m(\pi_{\mathbf{h}}\cap\aff(-\mathbf{b},-\mathbf{b}+\boldsymbol{\omega}_i))&=
\bigl\{-\mathbf{b}+\lambda\boldsymbol{\omega}_i\mid
\langle(A^\ast)^{-m}\mathbf{h},-\mathbf{b}+\lambda\boldsymbol{\omega}_i\rangle
=1\biggr\}
\\
\label{eq6}
&=\biggl\{\boldsymbol{\omega}_i
\frac{1+\langle(A^\ast)^{-m}\mathbf{h},\mathbf{b}\rangle}
{\langle(A^\ast)^{-m}\mathbf{h},\boldsymbol{\omega}_i\rangle}\biggr\}.
\end{align}
Hence the inequality~\eqref{eq5} is equivalent to the inequality
$$
\inf_{m\in\mathbb{Z}}\,\prod_{i=1}^n
\frac{1+\langle(A^\ast)^{-m}\mathbf{h},\mathbf{b}\rangle}
{\langle(A^\ast)^m\mathbf{h},\boldsymbol{\omega}_i\rangle}>0,
$$
which, in virtue of~\eqref{eq4}, implies that
\begin{equation}
\label{eq7}
\sup_{m\in\mathbb{Z}}\,\prod_{i=1}^n\langle(A^\ast)^m\mathbf{h},
\boldsymbol{\omega}_i\rangle<\infty.
\end{equation}

Let us assume that~$\mathbf{e}_i$ corresponds to~$\lambda_i$ and that
$|\lambda_1|\leq\nobreak\dots\leq\nobreak|\lambda_n|$. We denote by $\mu_1<\dots<\mu_k$, $k\leq n$,
the set of absolute values of $\lambda_1,\nobreak\dots,\nobreak\lambda_n$, and for every
$i\in\{1,\dots,k\}$ we denote by~$M_i$ the invariant subspace of~$A$ corresponding to~$\mu_i$. We
have a graduation $\mathbb{R}^n=M_1\oplus\dots\oplus M_k$.

Let us prove the following two statements:

(i) $M_i$ contains exactly $\dim M_i$ edges of the cone~$\mathscr{C}$;

(ii) the matrix of the operator~$A$ does not have any nontrivial Jordan cells.

Set
$$
L_i=\bigoplus_{j=1}^iM_j,\qquad
L_i'=\bigoplus_{j=i}^kM_j.
$$
We have filtrations $\mathbb{R}^n=L_k\supset L_{k-1}\supset\dots\supset L_1=M_1$ and
$\mathbb{R}^n=L_1'\supset L_2'\supset\nobreak\dots\supset\nobreak L_k'=M_k$. Set
$r_i=\min\{r\colon\boldsymbol{\omega}_i\in L_r\}$. Then, by Lemma~\ref{lem4}, for every
$i\in\{1,\dots,n\}$ there is an $l_i\in\mathbb{Z}$, $l_i\geq 0$, such that
\begin{equation}
\label{eq8}
\langle (A^\ast)^m\mathbf{h},\boldsymbol{\omega}_i\rangle\asymp
\mu_{r_i}^mm^{l_i},\qquad
m\to+\infty.
\end{equation}
Hence
\begin{equation}
\label{eq9}
\prod_{i=1}^n\langle(A^\ast)^m\mathbf{h},\boldsymbol{\omega}_i\rangle\asymp
\biggl(\,\prod_{i=1}^n\mu_{r_i}^m\biggr)\biggl(\,\prod_{i=1}^nm^{l_i}\biggr)=
\biggl(\,\prod_{i=1}^n\frac{\mu_{r_i}}{|\lambda_i|}\biggr)^m
\biggl(\,\prod_{i=1}^nm^{l_i}\biggr),\quad
m\to+\infty.
\end{equation}

Without loss of generality we may assume that $r_1\leq r_2\leq\dots\leq r_n$. Then $i\leq\dim
L_{r_i}$, which implies that $\mu_{r_i}\geq|\lambda_i|$, where the equality holds only if $i>\dim
L_{r_i-1}$. Thus, it follows from~\eqref{eq7} and~\eqref{eq9} that for every $i\in\{1,\dots,n\}$
$$
\dim L_{r_i-1}<i\leq\dim L_{r_i}.
$$
Hence for every $j\in\{1,\dots,k\}$ and every $i\leq\dim L_j$
$$
\dim L_{r_i-1}<i\leq\dim L_j,
$$
that is, $r_i\leq j$. Consequently, if $i\leq\dim L_j$, then $\boldsymbol{\omega}_i\in L_j$. This
means that~$L_j$ contains exactly $\dim L_j$ vectors~$\boldsymbol{\omega}_i$. Therefore,
if~$A^\ast$ has a nontrivial Jordan cell, then there is an~$i$ such that $l_i>0$ in~\eqref{eq8},
which, in virtue of~\eqref{eq9}, contradicts~\eqref{eq7}. Hence neither~$A$, nor~$A^\ast$ have
nontrivial Jordan cells, and~(ii) is proved.

Set $r_i'=\min\{r\colon\boldsymbol{\omega}_i\in L_r'\}$. Then
$$
\langle (A^\ast)^{-m}\mathbf{h},\boldsymbol{\omega}_i\rangle
\asymp\mu_{r_i'}^{-m},
\qquad
m\to+\infty.
$$
Renumbering, if necessary, the $\boldsymbol{\omega}_i$, we assume that $r_1'\leq r_2'\leq\dots\leq
r_n'$. Arguments, similar to those for~$L_j$, show that each~$L_j'$ contains exactly $\dim L_j'$
vectors~$\boldsymbol{\omega}_i$. But $L_j\cap L_j'=M_j$ and $\dim L_j+\dim L_j'=n+\dim M_j$, hence
$M_j$ contains exactly $\dim M_j$ vectors~$\boldsymbol{\omega}_i$. This proves~(i).

If $\mathbf{b}\neq\boldsymbol{0}$, then due to~(P1),~(P2), and~\eqref{eq4} the point~$\mathbf{b}$
is contained in some~$M_{j_0}$, corresponding to exactly two complex conjugate eigenvalues with
their absolute values equal to~$1$. At the same time, by~(i), there are
some~$\boldsymbol{\omega}_{i_1}$ and~$\boldsymbol{\omega}_{i_2}$ contained in~$M_{j_0}$, that
is,~$M_{j_0}$ contains a two-dimensional face~$\mathscr{F}$ of~$\mathscr{C}$. But since
$\langle(A^\ast)^m\mathbf{h},\mathbf{v}_{-m}\rangle= \langle\mathbf{h},\mathbf{v}_0\rangle$ for all
$m\in\mathbb{Z}$, we have either $\langle(A^\ast)^m\mathbf{h},\mathbf{x}\rangle\geq
\langle\mathbf{h},\mathbf{v}_0\rangle$ for all $\mathbf{x}\in K$, $m\in\mathbb{Z}$ and
$i\in\{1,\dots,n\}$, or $\langle(A^\ast)^m\mathbf{h},\mathbf{x}\rangle\leq
\langle\mathbf{h},\mathbf{v}_0\rangle$ for all $\mathbf{x}\in K$, $m\in\mathbb{Z}$ and
$i\in\{1,\dots,n\}$. Therefore, all the quantities
$\langle(A^\ast)^m\mathbf{h},\boldsymbol{\omega}_i\rangle$ for all $m\in\mathbb{Z}$ and
$i\in\{1,\dots,n\}$ have same signs. But this is not so, since the restriction of~$A$ to~$M_{j_0}$
is a rotation operator, distinct from the unity. Hence $\mathbf{b}=\boldsymbol{0}$.

Due to~(ii) and~(P2) each space~$M_j$ corresponds either to a totally real~$\lambda$, and then all
the vectors in~$M_j$ are eigenvectors, or to a pair of complex conjugate~$\lambda$
and~$\overline\lambda$. If~$\lambda$ is not a positive real number, then the arguments, similar to
those used for~$M_{j_0}$, lead to a contradiction with the fact that all the quantities
$\langle(A^\ast)^m\mathbf{h},\boldsymbol{\omega}_i\rangle$ for all $m\in\mathbb{Z}$ and
$i\in\{1,\dots,n\}$ have same signs. Therefore, all the eigenvalues of~$A$ are real and positive,
and~$\boldsymbol{\omega}_i$ are eigenvectors of $A$. Hence $A(\mathscr{C})=\mathscr{C}$. If the
characteristic polynomial of~$A$ is reducible, then~$A$ has a nontrivial invariant integer
subspace. In case all the eigenvalues of~$A$ are pairwise distinct, this contradicts the
irrationality of~$\mathscr{C}$. Thus, if the eigenvalues of~$A$ are pairwise distinct, then $A$ is
a hyperbolic operator.
\end{nproof}

\section{Three-dimensional case}
\label{sec6}

In this section we explain how to improve Theorems~\ref{th2} and~\ref{th3} in the three-dimensional
case and obtain Theorem~\ref{th5}.

The implication 1)$\Rightarrow$2) in Theorem~\ref{th5}, same as in Theorems~\ref{th1}
and~\ref{th4}, is obvious.

Suppose that the statement~2) of Theorem~\ref{th5} holds. Then, by Theorem~\ref{th2}, there is an
$\mathbf{A}\in\Aff_n(\mathbb{Z})$ establishing a nontrivial shift of
$\bigcup_{i\in\mathbb{Z}}\oSt_{\mathbf{v}_i}$ along itself. To prove~1) it suffices to show that
$\mathbf{A}\in\mathfrak{A}_0$ and that the eigenvalues of the linear component of~$\mathbf{A}$ are
pairwise distinct. Then we can and apply Theorem~\ref{th3}.

Let $\mathbf{A}\colon\mathbf{x}\mapsto A(\mathbf{x})+\mathbf{a}$, $A\in\SL_n(\mathbb{Z})$,
$\mathbf{a}\in\mathbb{Z}^n$. Suppose that~$A$ does not satisfy at least one of the statements~(P1)
and~(P2) or has two equal eigenvalues. It is clear that in this case all the eigenvalues of~$A$
have absolute values equal to~$1$, i.e. the matrix of~$A$ in the Jordan basis has one of the
following forms:
\begin{alignat*}{2}
1)\quad A&=\begin{pmatrix}
1 & 0 & 0 \\
0 & 1 & 0 \\
0 & 0 & 1
\end{pmatrix};&\qquad
5)\quad A&=\begin{pmatrix}
1 & 0 & 0 \\
0 & \lambda & 0 \\
0 & 0 & \lambda^{-1}
\end{pmatrix};
\\
2)\quad A&=\begin{pmatrix}
1 & 0 & 0 \\
0 & 1 & 1 \\
0 & 0 & 1
\end{pmatrix};&\qquad
6)\quad A&=\begin{pmatrix}
1 & 0 & 0 \\
0 & \cos\varphi & \sin\varphi \\
0 & -\sin\varphi & \cos\varphi
\end{pmatrix};
\\
3)\quad A&=\begin{pmatrix}
1 & 1 & 0 \\
0 & 1 & 1 \\
0 & 0 & 1
\end{pmatrix};&\qquad
7)\quad A&=\begin{pmatrix}
-1 & 0 & 0 \\
0 & \cos\varphi & \sin\varphi \\
0 & -\sin\varphi & \cos\varphi
\end{pmatrix}.
\\
4)\quad A&=\begin{pmatrix}
1 & 0 & 0 \\
0 & -1 & 1 \\
0 & 0 & -1
\end{pmatrix};
\end{alignat*}

Clearly, $\mathbf{v}_m=\mathbf{A}^m(\mathbf{v}_0)=
A^m(\mathbf{v}_0)+(A^{m-1}+\dots+A+E)\mathbf{a}$. In each case this relation allows to write down
the asymptotic of coordinates of~$\mathbf{v}_m$ and come to a contradiction either with the cone's
irrationality, or with the fact that all the~$\mathbf{v}_i$ lie in the interior of~$\mathscr{C}$
and tend to its boundary as $i\to\pm\infty$. The only difficulty of this argument is the number of
possibilities, so we skip the details.

\section{A relation to the Littlewood and Oppenheim conjectures}
\label{sec7}

The following two conjectures are classical.

\begin{littlewood}
If $\alpha,\beta\in\mathbb{R}$, then $\inf_{m\in\mathbb{N}}m\|m\alpha\|\,\|m\beta\|=0$,
where~$\|\,\cdot\,\|$ denotes the distance to the nearest integer.
\end{littlewood}

\begin{oppenheim}
If $n\geq 3$ and $\langle\mathbf{L}_1,\,\cdot\,\rangle,\dots, \langle\mathbf{L}_n,\,\cdot\,\rangle$
are $n$~linearly independent linear forms on~$\mathbb{R}^n$ such that
$$
\inf_{\mathbf{x}\in\mathbb{Z}^n\setminus\{\boldsymbol{0}\}}
|\langle\mathbf{L}_1,\mathbf{x}\rangle\dots
\langle\mathbf{L}_n,\mathbf{x}\rangle|>0,
$$
then the lattice $\bigl\{(\langle\mathbf{L}_1,\mathbf{x}\rangle,\dots,
\langle\mathbf{L}_n,\mathbf{x}\rangle)\mid \mathbf{x}\in\mathbb{Z}^n \bigr\}$ is algebraic (that
is, similar modulo the action of the group of diagonal $(n\times n)$-matrices to the lattice of a
complete module of a totally real algebraic number field of degree~$n$).
\end{oppenheim}

As is known (see~\cite{11}), the three-dimensional Oppenheim conjecture implies the Littlewood
conjecture. In~\cite{12} and~\cite{13} an attempt was made to prove the Oppenheim conjecture,
however, there was an essential gap in the proof. Thus, both conjectures remain unproved.

The results of the current paper together with those of~\cite{14} and~\cite{15} allow to
reformulate the Oppenheim conjecture in terms of Klein polyhedra. In~\cite{14} and~\cite{15} along
with the concept of determinant of an edge star (see Definition~\ref{def6}) the concept of
determinant of a face is considered.

\begin{defin}
\label{def8} Let $F$ be an arbitrary $(n-1)$-dimensional face of a sail $\Pi\subset\mathbb{R}^n$
and let $\mathbf{v}_1,\dots,\mathbf{v}_m$ be the vertices of~$F$. Then we define the {\it
determinant\/} of~$F$ as
$$
\det F= \sum_{1\leq i_1<\dots<i_n\leq m}
|\det(\mathbf{v}_{i_1},\dots,\mathbf{v}_{i_n})|\,. $$
\end{defin}

In these papers the following theorem is proved.

\begin{Theorem}
\label{th7} Let $\langle\mathbf{L}_1,\,\cdot\,\rangle,\dots,\langle \mathbf{L}_n,\,\cdot\,\rangle$
be $n$~linearly independent irrational linear forms on~$\mathbb{R}^n$ and let
$$
\mathscr{C}=\bigl\{\mathbf{x}\in\mathbb{R}^n\mid
\langle\mathbf{L}_i,\mathbf{x}\rangle\geq 0, \ i=1,\dots,n \bigr\}.
$$
Then the following statements are equivalent:

{\rm 1)} $\inf_{\mathbf{x}\in\mathbb{Z}^n\setminus\{\boldsymbol{0}\}}
|\langle\mathbf{L}_1,\mathbf{x}\rangle\dots \langle\mathbf{L}_n,\mathbf{x}\rangle|>0$;

{\rm 2)} the faces and the edge stars of the vertices of the sail~$\Pi$ generated by~$\mathbb{Z}^n$
and~$\mathscr{C}$ have uniformly bounded determinants.
\end{Theorem}

Due to Corollary~\ref{cor1}, the statement~2) of Theorem~\ref{th7} is equivalent to the fact that
there are finitely many affine types of complete stars of the sail's vertices (a {\it complete
star\/} of a vertex is the union of all the faces incident to this vertex).

Using Theorem~\ref{th4}, we get the following reformulation of the Oppenheim conjecture.

\begin{sailed_oppenheim}
Let $\langle\mathbf{L}_1,\,\cdot\,\rangle,\dots, \langle\mathbf{L}_n,\,\cdot\,\rangle$ be
$n$~linearly independent irrational linear forms on~$\mathbb{R}^n$, $n\geq 3$, and let
$$
\mathscr{C}=\bigl\{\mathbf{x}\in\mathbb{R}^n\mid
\langle\mathbf{L}_i,\mathbf{x}\rangle\geq 0, \ i=1,\dots,n\bigr\}.
$$
Suppose that the sail~$\Pi$ corresponding to~$\mathbb{Z}^n$ and~$\mathscr{C}$ has finitely many
affine types of complete stars of vertices. Then there is an unbounded in both directions (as a
subset of~$\mathbb{R}^n$) $\mathfrak{A}_1$-periodic chain $\{\mathbf{v}_i\}_{i\in\mathbb{Z}}$ of
vertices of~$\Pi$.
\end{sailed_oppenheim}

We remind that for $n=3$ (which is the most interesting case, since the three-dimensional Oppenheim
conjecture implies the Littlewood conjecture) one can replace $\mathfrak{A}_1$ by
$\Aff_3(\mathbb{Z})$. This allows to reformulate the Oppenheim conjecture for $n=3$ as follows.

First, we define a graph~$\mathscr{G}(\Pi)$. As the set of its vertices we take the set of pairs
$(F,\mathbf{v})$, where $F$ is a face of~$\Pi$ and $\mathbf{v}$ is a vertex of~$F$; we connect two
distinct vertices $(F,\mathbf{v})$ and $(G,\mathbf{w})$ of~$\mathscr{G}(\Pi)$ with an edge if,
firstly, $[\mathbf{v},\mathbf{w}]$ is a common edge of~$F$ and~$G$, and secondly, the bypass of~$F$
from~$\mathbf{v}$ to~$\mathbf{w}$ is counter-clockwise (with respect to the outer normal to the
Klein polyhedron). The graph~$\mathscr{G}(\Pi)$ is obviously planar and each of its vertices is
incident to exactly three edges.

Next, we define a colouring of~$\mathscr{G}(\Pi)$. In contrast to the colourings
of~$\mathscr{G}_k(\Pi)$ considered above, we colour both the vertices and the edges
of~$\mathscr{G}(\Pi)$. Let $(F,\mathbf{v})$ and $(G,\mathbf{w})$ be arbitrary vertices
of~$\mathscr{G}(\Pi)$. Let~$\mathbf{a}$ and $\mathbf{b}$ be the vertices of~$F$ next
to~$\mathbf{v}$, and let~$\mathbf{c}$ and $\mathbf{d}$ be the vertices of $G$ next to~$\mathbf{w}$.
Suppose that there is an $\Aff_3(\mathbb{Z})$-operator taking the union of the faces of~$\Pi$
incident to at least one of the vertices~$\mathbf{a}$, $\mathbf{v}$, and~$\mathbf{b}$ to the union
of the faces of~$\Pi$ incident to at least one of the vertices~$\mathbf{c}$,~$\mathbf{w}$,
and~$\mathbf{d}$. Then we colour $(F,\mathbf{v})$ and $(G,\mathbf{w})$ identically. With the edges
we do a similar thing. Let $((F_1,\mathbf{v}_1),(F_2,\mathbf{v}_2))$ and
$((G_1,\mathbf{w}_1),(G_2,\mathbf{w}_2))$ be arbitrary edges of~$\mathscr{G}(\Pi)$. Let
$\mathbf{a}_1$ be the vertex of $F_1$ next to~$\mathbf{v}_1$ and distinct from~$\mathbf{v}_2$, let
$\mathbf{a}_2$ be the vertex of $F_2$ next to~$\mathbf{v}_2$ and distinct from~$\mathbf{v}_1$, let
$\mathbf{b}_1$ be the vertex of~$G_1$ next to~$\mathbf{w}_1$ and distinct from~$\mathbf{w}_2$, and
let $\mathbf{b}_2$ be the vertex of~$G_2$ next to~$\mathbf{w}_2$ and distinct from~$\mathbf{w}_1$.
Suppose that there is an $\Aff_3(\mathbb{Z})$-operator taking the union of the faces of~$\Pi$
incident to at least one of the vertices~$\mathbf{a}_1$, $\mathbf{v}_1$, $\mathbf{v}_2$,
and~$\mathbf{a}_2$ to the union of the faces of~$\Pi$ incident to at least one of the
vertices~$\mathbf{b}_1$, $\mathbf{w}_1$, $\mathbf{w}_2$, and~$\mathbf{b}_2$. Then we colour
$((F_1,\mathbf{v}_1),(F_2,\mathbf{v}_2))$ and $((G_1,\mathbf{w}_1),(G_2,\mathbf{w}_2))$
identically.

Due to Theorem~\ref{th5}, the fact that there is a chain of vertices of~$\Pi$ with its images
in~$\mathscr{G}_2(\Pi)$ and~$\mathscr{G}_3(\Pi)$ having periodic $\Aff_3(\mathbb{Z})$-colourings is
equivalent to the fact that there is a chain of vertices of~$\mathscr{G}(\Pi)$ having periodic
``vertex-edge'' colouring. Hence for $n=3$ the Oppenheim conjecture is equivalent to the following
statement.

\begin{sailed_oppenheim_3dim}
Let $\langle\mathbf{L}_1,\,\cdot\,\rangle$, $\langle\mathbf{L}_2,\,\cdot\,\rangle$,
$\langle\mathbf{L}_3,\,\cdot\,\rangle$ be linearly independent irrational linear forms
on~$\mathbb{R}^3$ and let
$$
\mathscr{C}=\bigl\{\mathbf{x}\in\mathbb{R}^n\mid
\langle\mathbf{L}_i,\mathbf{x}\rangle\geq 0, \ i=1,2,3\bigr\}.
$$
Suppose that the colouring of the graph~$\mathscr{G}(\Pi)$ of the sail~$\Pi$ corresponding
to~$\mathbb{Z}^3$ and~$\mathscr{C}$ involves only a finite number of colours. Then there is an
unbounded in both directions (in the natural metric of~$\mathscr{G}(\Pi)$) periodically coloured
chain of vertices of~$\mathscr{G}(\Pi)$.
\end{sailed_oppenheim_3dim}

%\end{fulltext}

\end{document}